\input amstex.tex
\input amsppt.sty
\magnification1200
\hoffset4truemm
\hsize152truemm
\TagsOnRight
\catcode`\@=11
\mathsurround 1.6pt
\font\Bbf=cmbx12
\let\logo@\relax
\nopagenumbers
\title
Mixed hook-length formula
\endtitle
\author
Maxim Nazarov
\endauthor 


\def\ba{{\bar a}}
\def\bi{{\bar\imath}}
\def\bj{{\bar\jmath}}
\def\bx{{\boxed{\phantom{\square}}\kern-.4pt}}

\def\CC{{\Bbb C}}
\def\com{\ts,\hskip-.5pt}
\def\crc{{\raise.24ex\hbox{$\sssize\kern.1em\circ\kern.1em$}}}
\def\CS{{\CC\ts S}}

\def\de{\delta}
\def\dep{\de^{\ts\prime}}

\def\enddemos{{\ $\square$\enddemo}}

\def\Fc{F_\Lac\hskip-1pt}
\def\FL{F_\La}
\def\FM{\ts\overline{\!F}_\Mu}

\def\ga{\gamma}
\def\Ga{\Gamma}
\def\ge{\geqslant}

\def\la{\lambda}
\def\La{\Lambda}
\def\Lac{{\mathchoice
 {\textstyle\La^{\kern-1pt c}}{\textstyle\La^{\kern-2pt c}}
 {\scriptstyle\La^{\kern-2.5pt\hbox{$\scriptscriptstyle c$}}}{ch}}}
\def\Lap{{\mathchoice
 {\textstyle\La^{\kern-1pt\prime}}{\textstyle\La^{\kern-2pt\prime}}
 {\scriptstyle\La^{\kern-2.5pt\hbox{$\scriptscriptstyle\prime$}}}{ch}}}
\def\Lc{{\La\hskip-2.5pt\raise.5pt\hbox{$^c$}}}
\def\lc{{\ts,\hskip.95pt\ldots\ts,\ts\,}}
\def\le{\leqslant}

\def\Mu{{\operatorname{M}}}
\def\Mup{{\mathchoice
 {\textstyle\Mu^{\kern-1pt\prime}}{\textstyle\Mu^{\kern-1pt\prime}}
 {\scriptstyle\Mu^{\kern-1.5pt\hbox{$\scriptscriptstyle\prime$}}}{ch}}}
\def\mup{\mu^{\ts\prime}}

\def\ns{\hskip-1pt}
\def\nt{\noindent}

\def\om{\omega}
\def\ot{\otimes}

\def\si{\sigma}

\def\te{\theta}
\def\Te{\Theta}
\def\th{\tau}
\def\ts{\hskip1pt}

\def\vo{v_{\ts\Ga}}

\def\Wt{{W'}}

\def\ZZ{{\Bbb Z}}


\centerline{\Bbf Mixed hook-length formula for}
\smallskip
\centerline{\Bbf degenerate affine Hecke algebras}
\bigskip\medskip
\centerline{\bf Maxim Nazarov}
\bigskip\medskip

\abstract{
Take the degenerate affine Hecke algebra $H_{l+m}$ corresponding to the group
$GL_{l+m}$ over a $p$-adic field. Consider the $H_{l+m}$-module $W$ induced
from the tensor product of the evaluation modules over the algebras $H_l$
and $H_m$. The module $W$ depends on two partitions $\lambda$ of $l$ and $\mu$
of $m$, and on two complex numbers. There is a canonical operator
$J$ acting in $W$, it corresponds to the Yang $R$-matrix. 
The algebra $H_{l+m}$ contains the symmetric group algebra $\CS_{l+m}$
as a subalgebra, and $J$ commutes with the action of this subalgebra
in $W$. Under this action, $W$ decomposes into
irreducible subspaces according to the Littlewood-Richardson rule. We
compute the eigenvalues of $J$, corresponding to certain multiplicity-free
irreducible components of $W$. In particular, we give a formula for
the ratio of two eigenvalues of $J$, corresponding to the maximal and
minimal irreducible components. As an application of our results,
we derive the well-known hook-length formula for the dimension of the
irreducible $\CS_l$-module corresponding to $\la$.
}\endabstract


\nt
In this article we work with the degenerate affine Hecke algebra
$H_l$ corresponding to the general linear group $GL_l$ over a local
non-Archimedean field. This algebra was introduced by V.\,Drinfeld in
[D], see also [L]. The complex associative algebra $H_l$ is
generated by the symmetric group algebra $\CS_l$ and by the pairwise
commuting elements $x_1\lc x_l$ with the cross relations for $p=1\lc l-1$
and $q=1\lc l$
$$
\align
\si_{p,p+1}\ts x_q&=x_q\ts\si_{p,p+1}\ts,\quad q\neq p\ts,p+1\ts;
\\
\si_{p,p+1}\ts x_p&=x_{p+1}\ts\si_{p,p+1}-1\ts.
\endalign
$$
Here and in what follows $\si_{pq}\in S_l$ denotes
transposition of the numbers $p$ and $q$. 

For any partition $\la=(\la_1\ts,\la_2\ts,\ldots)$ of $l$ let $V_\la$ be
the corresponding irreducible $\CS_l$-\ts module.
There is a homomorphism $H_l\to\CS_l$ identical on the
subalgebra $\CS_l\subset H_l$ such that
$x_p\mapsto\si_{1p}+\ldots+\si_{p-1,p}$ for each
$p=1\lc l-1$. So $V_\la$ can be regarded as a $H_l$-\ts module.
For any number $z\in\CC$ there is also an automorphism of $H_l$
identical on the subalgebra $\CS_l\subset H_l$ such that
$x_p\mapsto x_p+z$ for each $p=1\lc l$. We will denote by 
$V_\la(z)$ the $H_l$-\ts module obtained by pulling $V_\la$ back through
this automorphism. The module $V_\la(z)$ is irreducible by definition.

Consider the algebra $H_l\ot H_m$ where both $l$ and $m$ are
positive integers. It is isomorphic to the subalgebra in $H_{l+m}$
generated by the transpositions $\si_{pq}$ where
$1\le p<q\le l$ or $l+1\le p<q\le l+m$, along with all the elements
$x_1\lc x_{l+m}$. For any partition $\mu$ of $m$ and any number $w\in\CC$
take the corresponding $H_m$-module $V_\mu(w)$.
Now consider the $H_{l+m}$-\ts module $W$ induced from the
$H_l\ot H_m$-\ts module $V_\la(z)\ot V_\mu(w)$.
Also consider the $H_{l+m}$-\ts module $\Wt$ induced from the
$H_m\ot H_l$-\ts module $V_\mu(w)\ot V_\la(z)$.
Suppose that $z-w\notin\ZZ$, then the modules $W$ and $\Wt$ are irreducible
and equivalent\ts; see [C]\ts. So there is a unique, up to scalar
multiplier, $H_{l+m}$-\ts intertwining operator $I:W\to \Wt$.
For a certain particular realization
of the modules $W$ and $\Wt$ we will give an explicit expression for
the operator $I$. In particular, this will fix the normalization of $I$.

Following [C]\ts, we will realize $W$ and $\Wt$ as certain left ideals
in the group algebra $\CS_{l+m}$. We will have $\Wt=W\ts\th$ where
$\th\in S_{l+m}$ is the permutation
$$
(\ts1\lc m\ts,m+1\lc m+l\ts)\,\ts\mapsto\,(\ts l+1\lc l+m\ts,1\lc l\ts)\,.
$$
Let $J:W\to W$ be the composition of the operator $I$ and
the operator $\Wt\to W$ of the right multiplication by $\th^{-1}$.
In \S2 we give an explicit expression for the operator $I$.
Using this expression makes the spectral analysis of the operator $J$
an arduous task\ts; cf.\ [AK]\ts.
However, our results are based on this expression.

The subalgebra $\CS_{l+m}\subset H_{l+m}$
acts in $W\subset\CS_{l+m}$ via left multiplication.
Under this
action the space $W$ splits into irreducible components
according to the Littlewood-Richardson rule [M]\ts.
The operator $J$ commutes with this action. Hence $J$
acts via multiplication by a certain complex number in every
irreducible component of $W$ appearing with multiplicity one.
In this article we compute these numbers for certain multiplicity-free
components of $W$, see Theorems 1 and 2.

For example,
there are two distinguished irreducible components of the  
$\CS_{l+m}$-module $W$
which always have multiplicity one. 
They correspond to the partitions
$$
\la+\mu=(\la_1+\mu_1\ts,\la_2+\mu_2\,,\ts\ldots)
\quad\text{and}\quad 
(\la'+\mup)\ts'=(\la'_1+\mup_1\ts,\la'_2+\mup_2\,,\ts\ldots)\ts'
\hskip-10pt
$$
where as usual $\la'=(\la'_1\ts,\la'_2\ts,\ldots)$ denotes the
conjugate partition. 
Denote by
$h_{\la\mu}(z,\hskip-1pt w)$
the ratio of the corresponding two eigenvalues of $J$, this ratio
does not depend on the normalization of the operator $I$.
Theorems 1 and 2 have the following
\proclaim\nofrills{Corollary\ts:}
$$
h_{\la\mu}(z,\hskip-1pt w)\ =\
\prod_{i,j}\,\frac
{\,z-w-\la'_j-\mu_i+i+j-1\ts}
{\ts z-w+\la_i+\mup_j-i-j+1\ts}
$$
where the product is taken over all
$i\ts,j=1\ts,2\ts,\ts\ldots$ such that $j\le\la_i\,,\mu_i$. 
\endproclaim

\nt
We derive this corollary in \S6.
Now identify partitions with
their Young diagrams. The condition $j\le\la_i\,,\mu_i$ means that
the box $(i\ts,j)$ belongs to the intersection of the diagrams
$\lambda$ and $\mu$. If $\la=\mu$ the numbers $\la_i+\lambda'_j-i-j+1$
are called [M] the {\it hook-lengths\/} of the diagram $\la$.
If $\la\neq\mu$ the numbers in the above fraction
$$
\la_i+\mup_j-i-j+1
\quad\text{and}\quad 
\la'_j+\mu_i-i-j+1
$$ 
may be called the {\it mixed hook-lengths\/}
of the first and of the second kind respectively.
Both these numbers are positive for any
box $(i\ts,j)$ in the intersection of $\la$ and $\mu$.

Let $h_\la$ denote the product of all $l$ hook-lengths of the Young diagram
$\la$. This product appears in the well-known {\it hook-length formula\/} [M]
for the dimension of the irreducible $\CS_l$-module $V_\la$:
$\dim V_\la={l\ts!}\ts/\ts{h_\la}$. In \S7 we derive
this formula from our Theorem 2. At the end of \S7 we discuss
applications of our Theorems 1 and 2 to the representaton theory
of affine Hecke algebras,~cf.~[LNT]\ts.

The present work arose from my conversations with B.\,Leclerc and
J.\ts-Y.\,Thibon. I am very grateful to them, and to the EPSRC for
supporting their visits to the University of York. 
I am also grateful to S.\,Kumar for a useful discussion.
I have been supported by the EC under the TMR grant FMRX-CT97-0100.


\medskip\nt{\bf\S1.}
Here we will collect several known facts about the irrdeducible
$\CS_l$-\ts modules.
Fix the chain $S_1\subset S_2\subset\cdots\subset S_l$
of subgroups with the standard embeddings.
There is a distinguished basis in the space $V_\la$ associated with this chain,
called the {\it Young basis\/}. Its vectors are labelled by the standard
tableaux [M] of shape $\la$. For every such a tableau $\La$ the basis vector
$v_\La\in V_\la$ is defined, up to a scalar multiplier, as follows.
For any $p=1\lc l-1$ take the tableau obtained from $\La$ by
removing each of the numbers $p+1\lc l$. Let the Young diagram $\pi$
be its shape. Then the vector $v_\La$ is contained in an irreducible
$\CS_p$-submodule of $V_\la$ corresponding to $\pi$.
Fix an $S_l$-invariant inner product $\langle\,\,,\,\rangle$ in $V_\la$.
The vectors $v_\La$ are then pairwise orthogonal. We will agree that
$\langle\ts v_\La,v_\La\ts\rangle=1$ for every $\La$.

There is an alternative definition [J] of the vector $v_\La\in V_\la$.
For each $p=1\lc l$ put $c_p=j-i$ if the number $p$ appears in the
$i$-\ts row and the $j$-th column of the tableau $\La$. The number
$c_p$ is called the {\it content} of the box of the diagram $\la$
occupied by $p$. Here on the left we show
the column tableau of shape $\la=(3\com2)$:

\vbox{\bigskip
$$
{\bx}
{\bx}
{\bx}
\phantom{\bx}
\phantom{\bx}
\phantom{\bx}
\phantom{\bx}
\phantom{\bx}
{\bx}
{\bx}
{\bx}
$$
\vglue-17.8pt
$$
{\bx}
{\bx}
\phantom{\bx}
\phantom{\bx}
\phantom{\bx}
\phantom{\bx}
\phantom{\bx}
\phantom{\bx}
{\bx}
{\bx}
\phantom{\bx}
$$
\vglue-44pt
$$
1
\kern9pt
3
\kern9pt
5
\kern81pt
0
\kern9pt
1
\kern9pt
2
$$
\vglue-18pt
$$
2
\kern9pt
4
\kern9pt
\phantom6
\kern78pt
\text{-1}
\kern9pt
0
\kern9pt
\phantom1
$$
\medskip\vskip3pt}

\noindent
On the right we indicated the contents of the boxes of the
Young diagram $\la=(3\com2)$. So here we have
$(\ts c_1\lc c_{\ts5}\ts)=(\ts0\ts,\text{-1}\ts,1\ts,0\ts,2\ts)$.
For each $p=1\lc l$ consider the image $\si_{1p}+\ldots+\si_{p-1,p}$
of the generator $x_p$ under the homomorphism $H_l\to\CS_l$. 

\proclaim\nofrills{Proposition 1\,({\rm\ts[J]\ts})\ts:}\
the element $v_\La$ of the $\CS_l$-\ts module $V_\la$ is determined,
\text{up to} a scalar multiplier, by equations
$(\ts\si_{1p}+\ns\ldots\ns+\si_{p-1,p}\ts)\cdot v_\La=c_p\ts v_\La$
for $p=1\lc l$.
\endproclaim

\nt
Take the diagonal matrix element of $V_\la$ corresponding
to the vector $v_\La$
$$
\FL\,=\,
\sum_{\si\in S_l}\,
\langle\,v_\La\com\ts\si\!\cdot\!\hskip.5pt v_\La\,\rangle
\,\si\,.
$$
\vskip-4pt\noindent
As a general property of matrix elements, we have the equality
$\FL^{\ts2}={l\ts!}/\dim V_\la\cdot\FL$.
There is an alternative expression for the element $\FL\in\CS_l$,
it goes back to \text{[C]\ts.} For any distinct $p\ts,q=1\lc l$
introduce the rational function of two complex variables $u\com v$
$$
f_{pq}(u\com v)=1-\si_{pq}\ts/(u-v)\,.
$$
These functions take values in the algebra $\CS_l$ and satisfy the relations
$$
\align
f_{pq}(u\com v)\,f_{pr}(u\com w)\,f_{qr}(v\com w)
&=
f_{qr}(v\com w)\,f_{pr}(u\com w)\,f_{pq}(u\com v)
\tag{1}
\\
f_{pq}(u\com v)\,f_{qp}(v\com u)
&=
1-(u-v)^{-2}\,.
\tag{2}
\endalign
$$
for all pairwise distinct indices $p\com q\com r$.
Introduce $l$ complex variables $z_1\lc z_l$.
Order lexicographically
the pairs $(p\com q)$ with $1\le p<q\le l$.
Define the rational function $\FL(z_1\lc z_l)$
as the ordered product of the functions
$f_{pq}(\ts z_p+c_p\com z_q+c_q\ts)$
over all the pairs $(p\com q)$.
Let $\Cal{Z}_\La$ be the vector subspace in $\CC^{\ts l}$ consisting of all
$l$-tuples $(z_1\lc z_l)$ such that $z_p=z_q$, whenever the numbers
$p$ and $q$ appear in the same row of the tableau $\La$.

\proclaim\nofrills{Proposition 2\,({\rm\ts[N]\ts})\ts:}\
the restriction of the rational function $\FL(z_1\lc z_l)$ to $\Cal{Z}_\La$
is regular at the origin ${(0\lc\ns0)}\in\CC^{\ts l}$,
and takes there the value $\FL$.
\endproclaim

\smallskip\nt{\bf\S2.}
Let us choose any standard tableau $\La$ of shape $\la$.
The $H_l$-\ts module $V_\la(z)$ can be realized as the left ideal in $\CS_l$
generated by the element $\FL$. The subalgebra $\CS_l\subset H_l$ acts here
via left multiplication. Due to
Proposition 1 and to the defining relations of $H_l$,
the action of the generators $x_1\lc x_l$ in this left ideal
can be then determined by setting
$x_p\cdot\FL=(c_p+z)\ts\FL$ for each $p=1\lc l$. 

Also fix any standard tableau $\Mu$ of shape $\mu$.
Let $d_q$ be the content of the box of the diagram $\mu$
occupied by $q=1\lc m$ in $\Mu$.
We will denote by $\overline{\!G}$ the image of any element $G\in\CS_m$
under the embedding $\CS_m\to\CS_{l+m}\ts:\,\si_{pq}\mapsto\si_{l+p,\ts l+q}$.
The $H_{l+m}$-module $W$ induced from the 
$H_l\ot H_m$-\ts module $V_\la(z)\ot V_\mu(w)$ can be realized
as the left ideal in $\CS_{l+m}$ generated by the product $\FL\ts\FM$.
The action of the generators $x_1\lc x_{l+m}$ in the latter left ideal
can be determined by setting
$$
\alignat2
x_p\cdot\FL\ts\FM&=\ts(c_p+z)\ts\FL\ts\FM
&&
\quad\text{for each}\quad p=1\lc l\,;
\tag{3}
\\
x_{l+q}\cdot\FL\ts\FM&=(d_q+w)\ts\FL\ts\FM
&&\quad\text{for each}\quad q=1\lc m\ts.
\endalignat
$$

Now introduce the ordered products in the symmetric group algebra $\CS_{l+m}$
$$
\align
R_{\La\Mu}(z\com w)\,\ts&=\!
\prod_{p=1,...,\ts l}^{\longrightarrow}
\biggl(\
\prod_{q=1,...,\ts m}^{\longleftarrow}
f_{p,l+q}(c_p+z\com d_q+w)
\biggr)\ts,
\\
R_{\La\Mu}^{\ts\prime}(z\com w)\,\ts&=\!
\prod_{p=1,...,\ts l}^{\longleftarrow}
\biggl(\
\prod_{q=1,...,\ts m}^{\longrightarrow}
f_{p,l+q}(c_p+z\com d_q+w)
\biggr)\ts;
\endalign
$$
the arrows indicate ordering of (non-commuting) factors. 
\!We keep to the \text{assumption} $z-w\notin\ZZ$. Applying Proposition 2
to the tableaux $\La,\Mu$ and using (1) repeatedly, we get
$$
\FL\ts\FM\ts R_{\La\Mu}(z\com w)=
R_{\La\Mu}^{\ts\prime}(z\com w)\ts\FL\ts\FM\,.
\tag{4}
$$
Hence the right multiplication in $\CS_{l+m}$ by
$R_{\La\Mu}(z\com w)$ preserves the left ideal~$W$\!.

The $H_{l+m}$-module $\Wt$ induced from the 
$H_m\ot H_l$-\ts module $V_\mu(w)\ot V_\la(z)$ can be then realized
as the left ideal in $\CS_{l+m}$ generated by the element
$\tau^{-1}\FL\ts\FM\ts\tau$,
so that $\Wt=W\tau$.
The action of the generators $x_1\lc x_{l+m}$ in $\Wt$
is determined~by
$$
\alignat2
x_q\cdot\tau^{-1}\FL\ts\FM\ts\tau&=\ts(d_q+w)\,\tau^{-1}\FL\ts\FM\ts\tau
&&
\quad\text{for each}\quad q=1\lc m\,;
\tag{5}
\\
x_{m+p}\cdot\tau^{-1}\FL\ts\FM\ts\tau&=(c_p+z)\,\tau^{-1}\FL\ts\FM\ts\tau
&&\quad\text{for each}\quad p=1\lc l\ts.
\tag{6}
\endalignat
$$
Consider the operator of right multiplication in $\CS_{l+m}$
by $R_{\La\Mu}(z\com w)\tau$.
Denote by $I$ the restriction of this
operator to the subspace $W\subset\CS_{l+m}$.
Due to (4) the image of the operator $I$ is contained in the subspace $\Wt$.

\proclaim\nofrills{Proposition 3\ts:}\,\
the operator $I:W\to \Wt$ commutes with the action of $H_{l+m}$.
\endproclaim

\demo{Proof}
The subalgebra $\CS_{l+m}\subset H_{l+m}$ acts in $W,\Wt$
via left multiplication\ts; the operator $I$ commutes with this action by 
definition. The left ideal $W$ is generated by the element $\FL\ts\FM$,
so it  suffices to check that
$x_p\cdot I\bigl(\FL\ts\FM\bigr)=I\bigl(\ts x_p\cdot\FL\ts\FM\bigr)$
for each $p=1\lc l+m$. Firstly consider the case $1\le p\le l$,
then by (3,4,5,6)
$$
\align
x_p\cdot I\bigl(\FL\ts\FM\bigr)
&=
\ts x_p\cdot\bigl(R_{\La\Mu}^{\ts\prime}(z\com w)\ts\tau\bigr)
\bigl(\tau^{-1}\FL\ts\FM\ts\tau\bigr)
\\
&=\bigl(R_{\La\Mu}^{\ts\prime}(z\com w)\ts\tau\bigr)
\bigl(\ts x_{\tau^{-1}(p)}\cdot\tau^{-1}\FL\ts\FM\ts\tau\bigr)
\\
&=\bigl(R_{\La\Mu}^{\ts\prime}(z\com w)\ts\tau\bigr)
(c_p+z)
\bigl(\tau^{-1}\FL\ts\FM\ts\tau\bigr)
=I\bigl(\ts x_p\cdot\FL\ts\FM\bigr)\,;
\endalign
$$
here we also used the defining relations of the algebra $H_{l+m}$.
For more details of this argument see [L]\ts.
The case $l+1\le p\le l+m$ can be considered similarly
\enddemos

\noindent
Consider
the operator of the right multiplication in $\CS_{l+m}$
by $R_{\La\Mu}(z\com w)$. This operator 
preserves the subspace $W$ due to (4). The restriction of this operator
to $W$ will be denoted by $J$.
The subalgebra $\CS_{l+m}\subset H_{l+m}$ acts in $W$
via left multiplication, and $J$ commutes with this action.
Regard $W$ as a $\CS_{l+m}$-\ts module only. Let $\nu$ be any partition of
$l+m$ such that the irreducible $\CS_{l+m}$-module $V_\nu$ appears in $W$
with multiplicity one. The operator $J:W\to W$ preserves the
subspace $V_\nu\subset W$ and acts there 
as multiplication by a certain number from $\CC$.
Denote this number by $r_\nu(z\com w)$; it depends on
$z$ and $w$ as a rational function of $z-w$,
and does not depend on the choice of the tableaux $\La$ and $\Mu$. Our aim
is to compute the eigenvalues $r_\nu(z\com w)$ of $J$ for certain $\nu$.

Before performing the computation, let us observe
one general property of the eigenvalues $r_\nu(z\com w)$. Similarly to the
defintion of the $H_{l+m}$-intertwining operator $I:W\to \Wt$,
one can define an operator $I':\Wt\to W$ as the restriction to $\Wt$
of the operator of the right multiplication in $\CS_{l+m}$
by $R_{\ts\Mu\La}(w\com z)\,\tau^{-1}$. The operator $I'$
commutes with the action of $H_{l+m}$ as well.
One can also consider the operator $J':\Wt\to\Wt$, defined
as the restriction to $\Wt$ of 
the operator of right multiplication in $\CS_{l+m}$
by $R_{\ts\Mu\La}(w\com z)$. There is a unique irreducible
$\CS_{l+m}$-submodule in $V_\nu^{\prime}\subset\Wt$
equivalent to $V_\nu\subset W$; actually here
we have $V_\nu^{\prime}=V_\nu\ts\tau$. Consider the
corresponding eigenvalue $r_\nu^{\,\prime}(w\com z)$ of the operator $J'$.

\proclaim\nofrills{Proposition 4\ts:}\,\
$$
r_\nu(z\com w)\,\ts r_\nu^{\,\prime}(w\com z)\,=\,
\prod_{i=1}^{\la'_1}\ \prod_{k=1}^{\mup_1}\ 
\frac{(z-w+\la_i-i+k)\,(z-w-\mu_k-i+k)}{(z-w+\la_i-\mu_k-i+k)\,(z-w-i+k)}\,.
$$
\endproclaim

\demo{Proof}
The product $r_\nu(z\com w)\,\ts r_\nu^{\,\prime}(w\com z)$ is
the eigenvalue of the composition $I'\!\crc\ns I:W\to W$,
corresponding to the subspace $V_\nu\subset W$. By definition,
this composition
is the operator of right multiplication in $W\subset\CS_{l+m}$ by the element
$$
\gather
R_{\La\Mu}(z\com w)\,\tau\,R_{\ts\Mu\La}(w\com z)\,\tau^{-1}\ =
\prod_{p=1,...,\ts l}^{\longrightarrow}
\biggl(\
\prod_{q=1,...,\ts m}^{\longleftarrow}
f_{p,l+q}(c_p+z\com d_q+w)
\biggr)\cdot\tau\ \times
\\
\prod_{q=1,...,\ts m}^{\longrightarrow}
\biggl(\
\prod_{p=1,...,\ts l}^{\longleftarrow}
f_{q,m+p}(d_q+w\com c_p+z)
\biggr)
\cdot\tau^{-1}\ =
\\
\prod_{p=1,...,\ts l}^{\longrightarrow}
\biggl(\
\prod_{q=1,...,\ts m}^{\longleftarrow}
f_{p,l+q}(c_p+z\com d_q+w)
\biggr)
\hskip2pt\cdot\hskip-6pt
\prod_{p=1,...,\ts l}^{\longleftarrow}
\biggl(\
\prod_{q=1,...,\ts m}^{\longrightarrow}
f_{l+q,p}(d_q+w\com c_p+z)
\biggr)
\\
=\ 
\prod_{p=1}^l\ 
\prod_{q=1}^m\ 
\Bigl(\ts1-(z-w+c_p-d_q)^{-2}\Bigr)\,.
\endgather
$$
The last equality has been obtained by using repeatedly the relations (2),
it shows that the composition $I'\!\crc\ns I$ is a scalar operator.
Now recall that the contents of the boxes in the same row of a Young diagram
increase by $1$, when moving from left to right. For the $i$-th row of
$\la$, the contents of the leftmost and rightmost boxes are $1-i$ and
$\la_i-i$ respectively. For the $k$-th row of $\mu$, the contents of
the leftmost and rightmost boxes are respectively $1-k$ and $\mu_k-k$.
Hence the right hand side of the last equality can be rewritten as
$$
\align
&
\prod_{p=1}^l\ 
\prod_{q=1}^m\ 
\biggl(\
\frac{z-w+c_p-d_q+1}{z-w+c_p-d_q}
\cdot
\frac{z-w+c_p-d_q-1}{z-w+c_p-d_q}
\,\biggr)\ =
\\
&
\prod_{i=1}^{\la'_1}\ \ts
\prod_{q=1}^m\ \ts 
\biggl(\
\frac{z-w+\la_i-i-d_q+1}{z-w+1-i-d_q}
\cdot
\frac{z-w-i-d_q}{z-w+\la_i-i-d_q}
\,\biggr)\ =
\\
&
\prod_{i=1}^{\la'_1}\ \ts
\prod_{k=1}^{\mup_1}\ \ts 
\biggl(\
\frac{z-w+\la_i-i+k}{z-w+\la_i-i-\mu_k+k}
\cdot
\frac{z-w-i-\mu_k+k}{z-w-i+k}
\,\biggr)
\quad\ \square
\endalign
$$
\enddemo

\smallskip\nt{\bf\S3.}
Choose any sequence $a_1\lc a_{\la'_1}\in\{1\ts,2\ts,\ts\ldots\}$
of pairwise distinct indices, we emphasize that 
this sequence needs not to be increasing. 
Here $\la'_1$ is the number of non-zero parts in the partition $\la$.
Consider the partition $\mu$ as an infinite sequence
with finitely many non-zero terms.
Define an infinite sequence $\ga=(\hskip1pt\ga_1\ts,\ga_2\ts,\,\ldots)$~by
$$
\gather
\ga_{a_i}=\mu_{a_i}+\la_i\,,\quad i=1\lc\la'_1\,;
\\
\ga_a=\mu_a\,,\quad a\neq a_1\lc a_{\la'_1}\,.
\endgather
$$
Suppose we have $\ga_1\ge\ga_2\ge\ldots\,$, so that 
$\ga$ is a partition of $l+m$.
Then the irreducible $\CS_{l+m}$-module $V_\ga$
corresponding to the partition $\ga$
\text{appears} in $W$ with multiplicity one.
Indeed, the multiplicity of $V_\ga$ in $W$ equals the multiplicity of
$V_\la$ in the $\CS_l$-module corresponding to the skew Young diagram
$\ga\ts/\mu$. The latter multiplicity is one by the
definition of $\ga$; see for instance [M].

We will evaluate the number $r_\ga(z\com w)$ by applying the operator $J$
to a particular vector in the subspace $V_\ga\subset W$.
Assume that $\La=\Lc$ is the column tableau
of shape $\la$; the tableau $\operatorname{M}$ will be still arbitrary.
The image of the action of the element $\Fc\FM$
in the irreducible $\CS_{l+m}$-\ts module $V_\ga$ is a
one-dimensional subspace. Let us describe this subspace explicitly.
The standard chain of subgroups
$S_1\subset S_2\subset\cdots\subset S_l$ corresponds to the natural ordering
of the numbers $1\ts,2\lc l$. Now consider the new chain of subgroups
$$
S_1\subset\ldots\subset S_m\subset S_{1+m}\subset\ldots\subset S_{l+m}
$$
corresponding to the ordering $l+1\lc l+m\,,1\lc l$. Notice that
the element $\FM$ belongs to the subgroup $S_m$ in this new chain.
Take the Young basis in the space $V_\ga$ associated with the new chain.
In particular, take the basis vector $\vo\in V_\ga$ corresponding to
the tableau $\Ga$ of shape $\ga$ defined as follows.
The numbers $l+1\lc l+m$
appear in $\Ga$ respectively in the same positions as
the numbers $1\lc m$ do in tableau $\Mu$.
Now for every positive integer $j$ consider all those parts of $\la$
which are equal to $j$. These are $\lambda_i$ with 
$i=\la'_{j+1}+1\,,\la'_{j+1}+2\ts\lc \la'_j$.
Let $f_1\lc f_{\la'_j-\la'_{j+1}}$ be the indices $a_i$ with $\la_i=j$,
arranged in the increasing order. By definition, the numbers appearing
in the rows $\la'_{j+1}+1\,,\la'_{j+1}+2\ts\lc \la'_j$ of the tableau $\Lc$,
will stand in the rows $f_1\lc f_{\la'_j-\la'_{j+1}}$ of $\Ga$ respectively.
For example, here for $\la=(3\ts,\ns2)$ and $\mu=(2\ts,\ns1)$ with
$a_1=2$ and $a_2=1$
we show a standard tableau $\Mu$ and the corresponding tableau $\Ga$:

\vbox{\bigskip
$$
{\bx}
{\bx}
\phantom{\bx}
\phantom{\bx}
\phantom{\bx}
\phantom{\bx}
\phantom{\bx}
{\bx}
{\bx}
{\bx}
{\bx}
$$
\vglue-17.8pt
$$
{\bx}
\phantom{\bx}
\phantom{\bx}
\phantom{\bx}
\phantom{\bx}
\phantom{\bx}
\phantom{\bx}
{\bx}
{\bx}
{\bx}
{\bx}
$$
\vglue-44pt
$$
1
\kern9pt
2
\kern81pt
6
\kern9pt
7
\kern9pt
2
\kern9pt
4
$$
\vglue-18pt
$$
3
\kern9pt
\phantom{0}
\kern81pt
8
\kern9pt
1
\kern9pt
3
\kern9pt
5
$$
\smallskip}

\proclaim\nofrills{Proposition 5\ts:}\ \,
with respect to the ordering
$l+1\lc l+m\,,1\lc l$
the tableau $\Ga$ is standard.
\endproclaim

\demo{Proof}
Reading the rows of the tableau $\Ga$ from left to right,
or reading its columns downwards,  
the numbers $l+1\lc l+m$ appear in the increasing order because
$\Mu$ is standard. These numbers will also appear before $1\lc l$.
Moreover, the numbers $1\lc l$
increase along each row of the tableau $\Ga$ by the definition of $\Lc$.
Now suppose that a column of $\Ga$ contains two
different numbers $p\ts,q\in\{1\lc l\ts\}$.
Let $a\ts,\ba$ be the corresponding rows\ts; assume that $a<\ba$. Then
$a=a_i$ and $\ba=a_\bi$ for certain indices $i\ts,\bi\in\{1\lc\la'_1\}$.
If $\la_i\ge\la_\bi$ then $p<q$ by
definition of $\Lc$.

Let $j,\bj$ be the columns corresponding to the numbers $p,q$ 
in the tableau $\Lc$. Suppose that $\la_i<\la_\bi$, then
$\mu_a>\mu_\ba$ because $\mu_a+\la_i\ge\mu_\ba+\la_\bi$. Since
$p$ and $q$ stand in the same column of the tableau $\Ga$,
we then have $j<\bj$ and $p<q$
\enddemos

\proclaim\nofrills{Proposition 6\ts:}\,\
the one-dimensional subspace $\Fc\FM\cdot V_\ga\subset V_\ga$
is obtained from the space $\CC\ts\vo$
by antisymmetrization relative to the columns of the tableau $\Lc$.
\endproclaim

\demo{Proof}
Let $S_\la$ be the subgroup in $S_l$ consisting of all permutations
which preserve the columns of the tableau $\Lc$ as sets. Let $Q\in\CS_l$
be the alternated sum of all elements from $S_\la$.
Put $V=\FM\cdot V_\ga$.
The subspace $V\subset V_\ga$ is spanned by the Young vectors, corresponding
to the tableaux which agree with $\Ga$ in the entries $l+1\lc l+m$.
The action of the element $\Fc$ in $V_\ga$ preserves the subspace $V$,
and the image $\Fc\cdot V\subset V$ is one-dimensional. Moreover, we have
$\Fc\cdot V=Q\cdot V$; see [JK].
It now remains to check that $Q\cdot\vo\neq0$.

By our choice of the tableau $\Ga$,
it suffices to consider the case when $\la$ consists 
of one column only. But then the element $Q\in\CS_l$ is central.
On the other hand, the vector $\vo\in V$ is $\CS_l$-\ts cyclic; see [C].
So $Q\cdot V\neq\{0\}$ implies $Q\cdot\vo\neq0$
\enddemos

\smallskip\nt{\bf\S4.}
Let $\ga$ be any of the partitions of $l+m$ described in the beginning of \S3.
Our main result is the following expression for the corresponding eigenvalue
$r_\ga(z\com w)$ of the operator $J:W\to W$.
This expression will be obtained by applying $J$ to the vector
$\Fc\FM\cdot\vo$ in $V_\ga\subset W$ and using Proposition 6.

\proclaim\nofrills{Theorem 1\ts:}
$$
r_\ga(z,\ns w)\ =\
\prod_{(i,j)}\,\frac
{\,z-w-\la'_j-\mu_{a_i}+a_i+j-1\,}
{\ts z-w-i+j\ts}
\hskip-30pt
$$
where the product is taken over all boxes $(i\com j)$ of the Young
diagram $\la$.
\endproclaim

\demo{Proof}
Using (4) and applying Proposition 2.12 of [N] to the tableau $\Mu$, we obtain
the equalities in the algebra $\CS_{l+m}$
$$
\Fc\FM\ts R_{\Lac\!\Mu}(z\com w)=\dim V_\la\ts/\ts{l\ts!}\ts\cdot
\Fc\FM\ts R_{\Lac\!\Mu}(z\com w)\ts\Fc=\dim V_\la\ts/\ts{l\ts!}\ \times
$$
\vskip-15pt
$$
\Fc\FM\,
\biggl(\ 
\prod_{p=1,...,\ts l}^{\longrightarrow}\ts
\frac
{\si_{l+1,p}+\ns\ldots\ns+\si_{l+m,p}-\ns c_p\ns-\ns z\ns+\ns w}
{-c_p\ns-\ns z\ns+\ns w}\,
\biggr)\ts
\Fc=\dim V_\la\ts/\ts{l\ts!}\ \times
$$
$$
\Fc\FM\,\ts
\biggl(\ 
\prod_{p=1,...,\ts l}\ts
\frac
{\si_{l+1,p}+\ns\ldots\ns+\si_{l+m,p}
+
\si_{1p}+\ns\ldots\ns+\si_{p-1,p}-\ns2c_p\ns-\ns z\ns+\ns w}
{-c_p\ns-\ns z\ns+\ns w}\,
\biggr)\,
\Fc\,;
$$
we have also used Proposition 1 of the present article, cf.\ [O]\ts.
Here in the last line the factors corresponding to $p=1\lc l$
pairwise commute.
By the same proposition applied to the partition $\ga$ instead of $\la$,
any Young vector in the $\CS_{l+m}$-\ts module $V_\ga$ is an eigenvector
for the action of the elements
$$
\si_{l+1,p}+\ns\ldots\ns+\si_{l+m,p}+\si_{1p}+\ns\ldots\ns+\si_{p-1,p}\ts\,;
\qquad
p=1\lc l\,.
$$
The vector $Q\cdot\vo\in V_\ga$ is a linear combination of the Young
vectors, corresponding to standard tableaux obtained from $\Ga$
by permutations from the subgroup $S_\la\subset S_l$.
The last expression for $\Fc\FM\ts R_{\Lac\!\Mu}(z\com w)$ now shows,
in particular, that the number $r_\ga(z,\ns w)$ factorizes with
respect to the columns of the Young diagram $\la$.

Firstly suppose that $\la$ consists of one column only.
Then the number $r_\ga(z,\ns w)$ is
easy to evaluate\ts; cf.\ [NT]\ts. Here we have $c_i=1-i$ for each
$i=1\lc l$.
Using the chain of subgroups $S_1\subset S_2\subset\cdots\subset S_l$
corresponding to ordering $l\lc\ns1$ we then get
$$
\Fc\,\cdot\! 
\prod_{i=1,...,\ts l}^{\longrightarrow}\ts
\bigl(\,
\si_{l+1,i}+\ns\ldots\ns+\si_{l+m,i}-\!1+i+u
\ts\bigr)\ts=
$$
$$
\Fc\,\cdot\! 
\prod_{i=1,...,\ts l}\ts
\bigl(\,
\si_{l+1,i}+\ns\ldots\ns+\si_{l+m,i}+\si_{li}+\ns\ldots\ns+\si_{i+1,i}-\!1+l+u
\ts\bigr)\,.
$$
Here in the last line the factors corresponding to $i=1\lc l$
pairwise commute. Their product commutes with any element
from the subalgebra $\CS_l\subset\CS_{l+m}$, and acts on the vector
$\vo\in V_\ga$ as multiplication by the number
$$
\prod_{i=1}^l\ 
\bigl(\,
\mu_{a_i}\!-a_i+l+u
\ts\bigr)\,.
\tag{7}
$$

Let us now apply this result to the $j$-th column of a general Young diagram
$\la$, consecutively for $j=1\lc\la_1$. For the general $\la$ the content of
the box $(i\com j)$ is $j\ns-\ns i$. According to our last expression for
$\Fc\FM\ts R_{\Lac\!\Mu}(z\com w)$ in the product (7) we then have to replace
$l\ts,\mu_{a_i}\ts,u$ by 
$\la^{\ts\prime}_j\ts,\mu_{a_i}+j-1\ts,2\ns-\ns2j\ns-\ns z\ns+\ns w$
respectively. Hence
$$
r_\ga(z,\ns w)\,=\,
\prod_{j=1}^{\la_1}\ \prod_{i=1}^{\la^{\ts\prime}_j}\   
\frac{\ts\mu_{a_i}+j-1-a_i+\la^{\ts\prime}_j+2\ns-\ns2j\ns-\ns z\ns+\ns w\ts}
{\ts i-j-z+w\ts}
\quad\ \square
$$
\vskip-5pt
\enddemo

\smallskip\nt{\bf\S5.}
Choose any sequence $b_1\lc b_{\la_1}\in\{1\ts,2\ts,\ts\ldots\}$
of pairwise distinct indices. Again, this sequence needs not to be
increasing. 
Let us now regard the partition $\mup$ conjugate to $\mu$ as an infinite
sequence with finitely many parts. Determine an infinite sequence
$\dep=(\hskip1pt\dep_1\ts,\dep_2\ts,\,\ldots)$ by
$$
\gather
\dep_{b_j}=\mup_{b_j}+\la'_j\,;\quad j=1\lc\la_1\,;
\\
\dep_b=\mup_b\,;\quad b\neq b_1\lc b_{\la_1}\,.
\endgather
$$
Suppose $\dep_1\ge\dep_2\ge\ldots\,$, so that 
$\dep$ is a partition of $l+m$. Define $\de$ as the partition 
conjugate to $\dep$.
The irreducible $\CS_{l+m}$-module $V_\de$
appears in $W$ with multiplicity one.
Take the corresponding eigenvalue $r_\de(z\com w)$ of the
\text{operator $J:W\to W$.}

\proclaim\nofrills{Theorem 2\ts:}
$$
r_\de(z\com w)\ =\
\prod_{(i,j)}\,\frac
{\,z-w+\la_i+\mup_{b_j}-i-b_j+1\ts}
{\ts z-w-i+j\ts}
\hskip-30pt
$$
where the product is taken over all boxes $(i\com j)$ of the Young
diagram $\la$.
\endproclaim

\demo{Proof}
Denote by $Z_\de$ the minimal central idempotent in
$\CS_{l+m}$ corresponding to the partition $\de$.  
Take the automorphism $\ast$ of the algebra $\CS_{l+m}$
such that $\si^\ast=\operatorname{sgn}(\si)\ts\si$; we have
$Z^{\ts\ast}_{\de}=Z^{\phantom{\ts\ast}}_{\dep}$ then.
Reflecting the tableaux $\La$ and $\Mu$ in their main diagonals
we get certain standard
tableaux of shapes $\la'$ and $\mup$ respectively; denote these
tableaux by $\La'$ and $\Mu^{\ts\prime}$. Then we have
$\FL^{\hskip.5pt\ast}=F_\Lap$ and
$F_\Mu^{\,\ast}=F_\Mup$.

On the other hand, by the definition of the number $r_\de(z\com w)$ we have
the equality
$$
Z_\de\,\FL\,\FM\,R_{\La\Mu}(z\com w)\ts=
\ts r_\de(z\com w)\,\ts Z_\de\,\FL\,\FM\,.
$$
By applying the automorphism $\ast$ to this equality we get
$$
Z_{\dep}\ts F_\Lap\ts\overline{\!F}_\Mup R^{\,\ast}_{\La\Mu}(z\com w)\ts=
r_\de(z\com w)\,\ts Z_{\dep}\ts F_\Lap\ts\overline{\!F}_\Mup.
$$
But $R^{\,\ast}_{\La\Mu}(z\com w)=R_{\Lap\!\Mup}\!(-z,\!-w)$ by definition.
Therefore by applying Theorem 1 to the partitions $\la',\mup$ instead of
$\la,\mu$ and choosing
$\ga=\dep$ we get
$$
r_\de(z,\ns w)\,=\,
\prod_{(i,j)}\,\frac
{\,z-w+\la_j+\mup_{b_i}-b_i-j+1\,}
{\ts z-w+i-j\ts}
\hskip-30pt
$$
where the product is taken over all boxes $(i\com j)$
of the diagram $\la'$. Equivalently, this product may be taken
over all boxes $(\ts j\ts,i)$ of the diagram $\la$
\enddemos

\smallskip\nt{\bf\S6.}
Let us now derive the Corollary stated in the beginning of this article.
We will use Theorems 1 and 2 in the simplest situation when
$a_i=i$ for every $i=1\lc\la'_1$ and $b_j=j$ for every $j=1\lc \la_1$.
Then we have
$\ga=\la+\mu$ and $\de=(\la'+\mup)\ts'$.
By Theorems 1 and 2,
$h_{\la\mu}(z\com w)=r_{\la+\mu}(z\com w)\ts/\ts r_{(\la'+\mup)\ts'}(z\com w)$
equals the product of the fractions
$$
\frac
{\,  z-w-\la'_j-\mu_i+i+j-1\ts}
{\ts z-w+\la_i+\mup_j-i-j+1\ts}
\tag{8}
$$
taken over all boxes $(i\com j)$ of the diagram $\la$. Consider those
boxes of $\la$ which do not belong to $\mu$.
These boxes form a skew Young diagram, let us denote it by $\om$.
To obtain the Corollary, it suffices to prove the following

\proclaim\nofrills{Proposition 7\ts:}\ \,
the product of the fractions $(8)\!$ over the boxes $(i\com j)$ of
\text{$\om\!$ equals $1\ts.$\!}
\endproclaim

\demo{Proof}
We will proceed by induction on the number of boxes in the diagram $\om$.
Let us write $u$ instead of $z-w$ for short.
When the diagram $\om$ is empty, the statement to prove is tautological.
Now let $(a,b)$ be any box of $\om$ such that by removing it from $\la$
we obtain again a Young diagram\ts;
then we have $\la_a=b$ and $\la'_b=a$.
By applying the induction hypothesis to
the last diagram instead of $\la$, we have to show that the product
$$
\gather
\prod_{j=\mu_a+1}^{b-1}
\frac
{\ts u\ns+\ns b\ns-\ns1\ns+\ns\mup_j\ns-\ns a\ns-\ns j\ns+\ns1}
{u\ns+\ns b\ns+\ns\mup_j\ns-\ns a\ns-\ns j\ns+\ns1}
\ \ \cdot
\prod_{i=\mup_b+1}^{a-1}
\frac
{\ts u\ns-\ns\mu_i\ns-\ns a+i+b-\ns1}
{\ts u\ns-\ns\mu_i\ns-\ns a+\ns1\ns+i+b-\ns1}
\ \ \times
\\
\times\ \
\frac
{\ts u\ns-\ns\mu_a\ns-\ns a+a+b-\ns1}
{\ts u\ns+\ns b\ns+\ns\mup_b-a-b+\ns1}
\tag{9}
\endgather
$$
equals $1$. Note that here we have $\mu_a<\la_a$ and $\mup_b<\la'_b$.

Suppose there is a box $(\ts\bi\ts,\bj\,)$ of $\om$ with
$\mu_a<\bj<\la_a$ and $\mup_b<\bi<\la'_b$\ts, such that by adding this box
to $\mu$ we obtain again a Young diagram. Then we have
$\mu_\bi=\bj-1$ and $\mup_\bj=\bi-1$.
For the last diagram instead
of $\mu$, the product (9) equals $1$ by the induction hypothesis. Then
it suffices to check the equality to~$1$~of
$$
\frac{\ts u+b-\ns1+\bi-\ns1-a-\bj+\ns1}{u+b+\bi-\ns1-a-\bj+\ns1}
\,\cdot\,
\frac{\ts u+b+\bi-a-\bj+\ns1}{u+b+\bi-a-\bj}
\ \ \times
$$
\vskip-5pt
$$
\frac{\ts u-\bj+\ns1-a+\bi+b-\ns1}{\ts u-\bj+\ns1-a+\ns1+\bi+b-\ns1}
\,\cdot\,
\frac{\ts u-\bj-a+\bi+b}{\ts u-\bj-a+\bi+b-\ns1}
\ \,.
$$
\vskip-10pt
$$
\text{But this product has form\ \ }
\frac{v-1}{v}\cdot\frac{v+1}{v}\cdot\frac{v}{v+1}\cdot\frac{v}{v-1}
\text{\ \ with\ \ }
v=u-a+b+\bi-\bj\ .
$$

\vskip-5pt
It remains to consider the case when there is no box $(\ts\bi\ts,\bj\,)$ in
$\om$ with the above listed
properties. Then $\mup_j=a-1$ for all $j=\mu_a+1\lc b-1$
and $\mu_i=b-1$ for all $i=\mup_b+1\lc a-1$.
Hence in this remaining case the product (9) equals
$$
\frac{\ts u+b-\ns1+a-\ns1-a-b+\ns1+\ns1}{\ts u+b+a-\ns1-a-\mu_a-\ns1+\ns1}
\,\cdot\,
\frac{\ts u-b+\ns1-a+\mup_b+\ns1+b-\ns1}{\ts u-b+\ns1-a+\ns1+a-\ns1+b-\ns1}
\ \ \times
$$
\vskip-5pt
$$
\times\ \
\frac
{\ts u-\mu_a\ns-\ns a+a+b-\ns1}
{\ts u\ns+\ns b\ns+\ns\mup_b-a-b+\ns1}
\ =\ 1\ \quad\square
$$
\enddemo

\smallskip\nt{\bf\S7.}
In this final section we derive from Theorem 2 the formula
$\dim V_\la={l\ts!}\ts/\ts{h_\la}$ for the dimension of the
irreducible $\CS_l$-module $V_\la$.
We will actually show that the coefficient ${l\ts!}/\dim V_\la$
in the relation $\FL^{\ts2}={l\ts!}/\dim V_\la\cdot\FL$ equals $h_\la$,
the product of the $l$ hook-lengths of the Young diagram $\la$.
We will use induction on the number of rows in $\la$.
If there is only one row in $\la$,
then $h_\la=l\ts!$ and $\dim V_\la=1$, so the desired equality is clear.
Let us now make the inductive assumption for $\la$, and consider
the Young diagram obtained by adding $m$ boxes to $\la$ in the row
$\la'_1+1$. Denote the new diagram by $\te$, we assume that $m\le\la_i$
for any $i=1\lc \la'_1$. Put $\mu=(m\com0\com0\com\,\ldots)$
and consider the eigenvalue $r_\te(z\com w)$ of the operator $J:W\to W$
corresponding to the multiplicity-free component $V_\te\subset W$.

In our case, there is only one standard tableau $\Mu$ of shape $\mu$.
Let $\Te$ be the unique standard tableau of shape $\te$
agreeing with $\La$ in the entries $1\lc l$; the numbers $l+1\lc l+m$ then
appear in the last row of $\Te$. By definition,
$$
F_\Te\cdot\FL\ts\FM\ts R_{\La\Mu}(z\com w)=
r_\te(z\com w)\cdot F_\Te\ts\FL\ts\FM=
h_\la\,m\ts!\ts\,r_\te(z\com w)\cdot F_\Te\,;
$$
the second equality here has been obtained using the inductive assumption.
On the other hand, due to Proposition~2 the matrix element
$F_\Te$ coincides with the value of the product
$\FL\ts\FM\ts R_{\La\Mu}(z\com w)$ at $z-w=\la'_1$. To make
the inductive step, it now remains to check that $h_{\ts\te}$ coincides
with the value of $h_\la\,m\ts!\ts\,r_\te(z\com w)$ at $z-w=\la'_1$.

Let us use Theorem 2 when $b_j=j$ for every $j=1\lc \la_1$.
With our particular choice of $\mu$, we then obtain that 
$r_\te(z\com w)$ equals the product over~$i=1\lc\la'_1$~of
$$
\prod_{j=1}^m\,
(z-w+\la_i-i-j+2)
\,\cdot\hskip-6pt
\prod_{j=m+1}^{\la_i}\hskip-4pt
(z-w+\la_i-i-j+1)
\,\cdot\ts
\prod_{j=1}^{\la_i}\,
\frac1{z-w-i+j}\ .
$$
Changing the running index $j$ to $\la_i-j+1$ in the products over $j=1\lc m$
and over $j=m+1\lc\la_i$ above, we obtain after cancellations the equality
$$
r_\te(z\com w)\,=\,
\prod_{i=1}^{\la'_1}\,
\frac{z-w+\la_i-i+1}{z-w+\la_i-m-i+1}\,.
$$
This equality shows that the
value $r_\te(z\com w)$ at $z-w=\la'_1$ coincides with the ratio
$$
\frac{h_{\ts\te}}{\ts h_\la\,m\ts!\ts}
\ =\ 
\prod_{i=1}^{\la'_1}\ \prod_{j=1}^m\ 
\frac{\la_i+\la'_1-i-j+2}{\la_i+\la'_1-i-j+1}
\ =\ 
\prod_{i=1}^{\la'_1}\ 
\frac
{\la'_1+\la_i-i+1}{\la'_1+\la_i-m-i+1}\
\quad\square
$$
\enddemo

\noindent
Let us make a few concluding remarks. Throughout \S\S1-6 we assumed 
that $z-w\notin\ZZ$\ts, then the $H_{l+m}$-\ts module $W$ is irreducible.
For $z-w\in\ZZ$,
our Corollary implies that the module $W$ is reducible if
$z-w$ is a mixed hook-length of the second kind relative to $\la$ and $\mu$,
and if $w-z$ is a mixed hook-length of the first kind.
When $\la=\mu$, our Corollary implies that the module $W$ is
reducible if $|\ts z-w\ts|$ is a hook-length of $\la$.
Moreover, then the module $W$ is irreducible [LNT] for all remaining values
$|\ts z-w\ts|$.

When $\la\neq\mu$ the $H_{l+m}$-\ts module $W$ maybe reducible while
neither $z-w$ is a mixed hook-length of the second kind, nor
$w-z$ is a mixed hook-length of the first kind. The irreducibility
criterion for the module $W$ with arbitrary $\lambda$ and $\mu$
has been also given in [LNT]. This work shows that the module 
$W$ is reducible if and only if the difference $z-w$ belongs to a certain 
finite subset $\Cal{S}_{\la\mu}\subset\ZZ$ determined in [LZ].
This subset satisfies the property $\Cal{S}_{\la\mu}=-\ts\Cal{S}_{\mu\la}$.

Denote by $\Cal{D}_{\la\mu}$ the union of the sets of all zeroes
and poles of the rational functions 
$r_{\la+\mu}(z\com w)\ts/\ts r_\nu(z\com w)$ in $z-w$, where $\nu$
ranges over all partitions $\ga$ and $\de$ described in \S3 and \S5
respectively. Then $\Cal{D}_{\la\mu}\subset\Cal{S}_{\la\mu}$.
Then $-\ts\Cal{D}_{\mu\la}\subset\Cal{S}_{\la\mu}$ also.
Using [LZ] one can demonstrate that if $\la'_1,\mup_1\le3$ then
$\Cal{D}_{\la\mu}\cup(-\ts\Cal{D}_{\mu\la})=\Cal{S}_{\la\mu}$.
However, 
$\Cal{D}_{\la\mu}\cup(-\ts\Cal{D}_{\mu\la})\neq\Cal{S}_{\la\mu}$
for general partitions $\la$ and $\mu$.
For example, if $\la=(8,\ns3,\ns2,\ns1,\ns0,\ns0,\ldots)$
and $\mu=(6,\ns4,\ns4,\ns0,\ns0,\ldots)$ then
$0\in\Cal{S}_{\la\mu}$
but $0\notin\Cal{D}_{\la\mu}\ts,\Cal{D}_{\mu\la}$.

For general $\la$ and $\mu$,
it would be interesting to point out for every $t\in\Cal{S}_{\la\mu}$
a partition $\nu$ of $l+m$,
such that the $\CS_{l+m}$-module $V_\nu$
appears in $W$ with multiplicity one, and such that
the ratio $r_{\la+\mu}(z\com w)\ts/\ts r_\nu(z\com w)$
has a zero or pole at $t$, as a rational function of $z-w$.


\bigskip\medskip
\centerline{\bf References}
\bigskip
\itemitem{[AK]}
{S.\,Alishauskas and P.\,Kulish},
{\it Spectral resolution of $SU(3)$-invariant solutions of the Yang\ts-Baxter
equation},
{J.\,Soviet Math.}
{\bf 35}
(1986),
2563--2574.

\itemitem{[C]}
{I.\,Cherednik},
{\it Special bases of irreducible representations of a degenerate
affine Hecke algebra},
{Funct.\,Anal.\,Appl.}
{\bf 20}
(1986),
76--78.

\itemitem{[D]}
{V\hskip-1pt.\,Drinfeld},
{\it\hskip-.3pt Degenerate affine Hecke algebras and Yangians},
{Funct.\,Anal.\,Appl.}
{\bf 20}
(1986),
56--58.

\itemitem{[JK]}  
{G.\,James and A.\,Kerber},
{\it The Representation Theory of the Symmetric Group},
Addison\ts-Wesley, Reading MA, 1981.

\itemitem{[J]}
{A.\,Jucys},
{\it Symmetric polynomials and the centre of the symmetric group ring\/},
{Rep.\,Math.\,Phys.}
{\bf 5}
(1974),
107--112.


\itemitem{[LNT]}
{B.\,Leclerc, M.\, Nazarov and J.\ts-Y.\,Thibon},
{\it Induced representations of affine Hecke algebras
and canonical bases of quantum groups},
{\tt math/0011074}.


\itemitem{[LZ]}
{B.\,Leclerc and A.\,Zelevinsky},
{\it Quasicommuting families of quantum Pl\"ucker coordinates},
{Amer.\,Math.\,Soc.\,Translat.}
{\bf 181}
(1998),
85--108.

\itemitem{[L]}  
G.\,Lusztig,
{\it Affine Hecke algebras and their graded version},
{J.\,Amer.\,Math.\,Soc.}
{\bf2} 
(1989),   
599--635.

\itemitem{[M]}   
{I. Macdonald},
{\it Symmetric Functions and Hall Polynomials},
Clarendon Press, Oxford, 1979.

\itemitem{[N]}
{M.\,Nazarov},
{\it Yangians and Capelli identities},
{Amer.\,Math.\,Soc.\,Translat.}
{\bf 181}
(1998),
139--163.

\itemitem{[NT]}
{M.\,Nazarov and V.\,Tarasov},
{\it On irreducibility of tensor products of Yangian modules},
{Intern.\,Math.\,Research Notices}
(1998),
125--150.

\itemitem{[O]}
{A.\,Okounkov},
{\it Young basis, \!Wick formula, and higher Capelli identities},
{\text{Intern.} Math.\,Research Notices}
(1996),
817--839.


\bigskip\medskip
\centerline
{Department of Mathematics, University of York, York YO1 5DD, England}
\centerline
{E\ts-mail\ts: {\tt mln1@york.ac.uk}}


\bye